\documentclass[12pt]{amsart}
\usepackage{amssymb,latexsym,amsmath}
\newtheorem{thm}{Theorem}[section]
\newtheorem{defn}[thm]{Definition}
\newtheorem{lemma}[thm]{Lemma}
\newtheorem{claim}[thm]{Claim}
\newtheorem{conj}[thm]{Conjecture}
\numberwithin{equation}{section}

\newcommand{\M}{\mathcal{M}_A}
\newcommand{\Le}{\mathcal{L}}
\newcommand{\ML}{\mathcal{L}(\mathcal{M},\psi)}
\newcommand{\R}{\mathbb{R}}
\newcommand{\N}{\mathbb{N}}
\newcommand{\q}{\mathbf{q}}
\newcommand{\ve}{\mathbf{x}}
\newcommand{\p}{\mathbf{p}}
\newcommand{\vex}{\mathbf{\tilde{x}}}
\newcommand{\Z}{\mathbb{Z}}
\newcommand{\un}{u_{\mathbf{x}}}

\newcommand{\ep}{\epsilon}

\newcommand{\ee}{\mathbf{e}}
\newcommand{\vv}{\mathbf{v}}
\newcommand\da{Diophantine approximation}
\newcommand\di{Diophantine}
\newcommand\kh{Khintchine}
\newcommand\hy{hyperplane}
\newcommand{\SL}{\operatorname{SL}}
\newcommand{\GL}{\operatorname{GL}}
\newcommand{\al}{\alpha}
\newcommand{\La}{\Lambda}
\newcommand{\Ga}{\Gamma}
\newcommand{\vw}{\mathbf{w}}
\begin{document}
\title
{ A Khintchine-type theorem  for hyperplanes}
\author{Anish Ghosh}

\begin{abstract} We obtain the convergence case of a
\kh\ type theorem for a large class of \hy s. Our approach to the
problem is from a dynamical viewpoint, and we modify a method due
to Kleinbock and Margulis to prove the result.
\end{abstract}
\maketitle
\section {Introduction}

The purpose of this paper is to prove the convergence case of a \kh\
type theorem for \hy s. Let us first introduce some notation, and
recall terminology from the theory of metric \da. Vectors will be
denoted in boldface, i.e $\ve = (x_1,\dots,x_{k}) \in \R^{k}$,
further we will denote the vector $(1,x_1,x_2,\dots,x_{k}) \in \R^{k
+ 1} $ by $\vex$.\\
Let $M_{n,1} $ denote the set of matrices with $n$ rows, $1$ column,
and real entries. For $v
> 0$, and $n \in \N$, define the set $\mathcal{W}_{v}(n,1) $ to be
the set of all matrices $A = (\alpha_{i})_{0 \leq i \leq n-1} \in
M_{n,1} $ for which there are infinitely many $q \in \Z $ such that
\[ max_{i}|p_i + \alpha_{i}q  | < |q|^{-v} ~ \] for some $\p \in \Z^{n}
$. We further define the sets :
\begin{defn}\label{defK}
\[ \mathcal{W}_{v}^{+}(n,1) = \bigcup_{u > v}\mathcal{W}_{u}(n,1)
\]
and \[\mathcal{W}_{v}^{-}(n,1)= \bigcap_{u <
v}\mathcal{W}_{u}(n,1).\]
\end{defn}
Note that from Definition \ref{defK}, it follows that, $v_1 \leq
v_2 \Rightarrow \mathcal{W}_{v_2}(n,1) \subset
\mathcal{W}_{v_1}(n,1)$.
\begin{defn}\label{defA}
For $A  \in M_{n,1}$, denote by $ \M$, the hyperplane
\[\M =
\{(\ve,\vex A)|~\ve \in \R^{n - 1} \}.
\]
\end{defn}

\begin{defn}
For a submanifold $\mathcal{M}$ of $\R^n$, we denote by $\ML$, the
set:\\ $\{\ve \in \mathcal{M} |~ |p + \ve \cdot \q| <
\psi(\|\q\|)$~~for infinitely many $\q \in \Z^n$, and some $p \in
\Z$\}.
\end{defn}
Here, $\|\q\| = $max$\{|q_1|,\dots,|q_n|\}$ is the height of the
vector $\q$. Estimating the (induced) measure and Hausdorff
dimension of $\ML$ has been a major theme in the theory of \da. We
now state the classical \kh\ theorem, one of the first results along
these lines.
\begin{thm}\label{thm:khin}
Let $\psi : \N \to \R$ be any function.
~Then\\
$|\Le(\R^n,\psi)| = 0 ~$whenever $~\sum_{k = 1}^{\infty} k^{n -
1}\psi(k) < \infty$.
\end{thm}

where $|~| $ denotes Lebesgue measure. The above theorem is referred
to as the convergence case of the \kh-Groshev theorem and turns out
to be an immediate consequence of the Borel-Cantelli lemma. Since
the Borel-Cantelli lemma works in only one direction, the
complementary divergence case is more difficult and needs deeper
ideas  from \da\ and ergodic theory. In fact, the divergence case
also needs restrictions on the function $\psi$ in the case $n = 1$.
Classical results in this theory are due to A.\kh, A.Groshev and
W.Schmidt, among others. For more details on the history of the
theorem and the classical convergence and divergence case,
 we direct the reader to  \cite{S}, \cite{H} and \cite{D1}.\\
We now turn to the theory of metric \da\ on manifolds. The idea is
to ask questions of \kh-Groshev type while restricting to
submanifolds of $\R^n$. This turns out to be more difficult, and
needs more sophisticated techniques even for the convergence case.
For classical results and techniques in the theory of metric
Diophantine approximation, we refer the reader to \cite{BD} and
\cite{H}.\\
Our approach to the problem is from a dynamical point of view. In
\cite{KM}, D.Kleinbock and G.A.Margulis have developed a technique
which relates Diophantine questions on manifolds to quantitative
non-divergence of unipotent trajectories on $\SL(n,\R)/\SL(n,\Z)$.
This allowed them to settle a long standing conjecture of
Sprindzhuk. Subsequently, in \cite{BKM}, the authors modify these
techniques to prove the convergence case of a \kh-type result for
non-degenerate submanifolds of $\R^n$. Further, in \cite{K},
D.Kleinbock shows that Diophantine properties of affine
subspaces are inherited by their non-degenerate submanifolds.\\
We will use a modification of this method to establish a
convergence \kh\ theorem for affine hyperplanes.\\
Let us mention an important special case of the \kh\ set-up. A
submanifold $\mathcal{M}$ of $ \R^n$ is called extremal if the set
$\Le(\mathcal{M},\psi_{\ep})$ has measure $0$ for every $\ep
> 0$ where $\psi_{\ep}(k)= k^{-(n + \ep)}$.
[Note : If a point is $\psi_{\ep}$-approximable for some $\ep >
0$, it is referred to as Very Well Approximable.] Criteria for
extremality of subspaces were examined in detail in \cite{K}.
There, the following necessary and sufficient condition
for the extremality of \hy\ was laid out:\\
\begin{thm}\label{defK2}\cite{K}
$\M$ is extremal if and only if
\[A \notin \mathcal{W}_{n}^{+}(n,1).  \]
\end{thm}
In \cite{BBKM}(\S $6.5$, problem $5$), the authors pose the problem
of finding criteria for affine subspaces to satisfy a \kh\ theorem
for convergence and divergence. This question was also part of a
more specific list of conjectures in \cite{K}(\S $6.3$). In this
paper, we provide
 a partial answer by providing sufficient conditions
under which a \hy\ satisfies the convergence part of Khintchine's
theorem. Specifically, we will consider \hy s $\M$ where the matrix
$A$ satisfies the condition :
\begin{equation}
\label{defA2} A \notin \mathcal{W}_{n}^{-}(n,1).
\end{equation}

Note that $ \mathcal{W}_{n}^{+}(n,1) \subseteq
\mathcal{W}_{n}^{-}(n,1)$. Consequently, we are looking at the
complement of a larger set of matrices. Let us now state the main
result of this paper.
\begin{thm}\label{thmA3}
Let $\psi : \N \to \R$ be a non-increasing function, and assume that
$A \notin \mathcal{W}_{n}^{-}(n,1)$ . Then,
\[|\Le(\M,\psi)| = 0 \]~whenever
 $\sum_{k = 1}^{\infty}~ k^{n - 1}\psi(k) < \infty$.
\end{thm}

where $|~| $ denotes the volume measure on $M$.\\
Thus, in terminology borrowed from \cite{BD} if $A$ satisfies
(\ref{defA2}), $\M $ is a Groshev-type manifold for
convergence.\\
It is instructive to compare the condition in Theorem \ref{defK2},
with that in (\ref{defA2}). From Definition \ref{defK} and the
remark immediately after, it follows that if $ A $ does not
satisfy (\ref{defA2}), then $\M $ belongs to the null set of
non-extremal \hy s. The Hausdorff dimension of this set has been
computed in \cite{D2}, and is equal to $1$. Thus,
Theorem \ref{thmA3} applies to a very large class of \hy s.\\
Note that the set $\ML$ could be quite large, even though it has
measure zero. In case we take $\psi$ to be a decreasing function,
Dickinson and Dodson have proved the following theorem :
\begin{thm}\cite{DD}
Let $M$ be an $m$-dimensional extremal $C^{1}$ manifold embedded in
$ \R^n$. Let $\psi : \N \to \R$ be decreasing with the lower order
of $\frac{1}{\psi}$ denoted by $\lambda(\psi)$. Then for
$\lambda(\psi) \geq n$,
\[dim \ML \geq m - 1 + \frac{n + 1}{\lambda(\psi) + 1}\]
\end{thm}
where dim denotes Hausdorff dimension and the lower order of $\psi$
is defined as $\lambda(\psi) = \liminf_{n \to
\infty}\frac{\log\psi(n)}{\log(n)}$.\\  For example, if we choose a
$A$ as in (\ref{defA2}) and for $k \geq 2$, let $\psi(k) =
\frac{1}{k^n (\log k)^{1 + \epsilon}}$, then Theorem \ref{thmA3}
combined with the above result tells us that dim$(\Le(\M,\psi))
= n - 1$.\\
We do not examine the complementary divergence case here, although
we do think it should be true. More precisely, we make the following
conjecture:
\begin{conj}
Let $ A $ satisfy (\ref{defA2}), and $\psi : \N \to
\R $ be a non-increasing function. Then,\\
$\Le(\M,\psi) $ has full measure whenever $\sum_{k =
1}^{\infty}k^{n - 1}\psi(k) = \infty $.
\end{conj}

The divergence case has been proved in \cite{BBDD} for straight
lines passing through the origin.\\
The paper is organized as follows. In Section $2$, we state
definitions and reduce the result to two Theorems. Section's $3$
and $4$ deal with these Theorems and complete the proof.

\section {Reduction to Borel-Cantelli }
In what follows, $|~|$ will denote both the Lebesgue measure of a
subset of $\R^n$, as well as the
absolute value of a real number. The context will hopefully make the usage clear.\\
Fix an open ball $B \subseteq \R^{n - 1}$, a function $\psi : \N \to
\R^{+}$ and  a matrix $A =
(\alpha_0,\alpha_1,\dots,\alpha_{n-1})^{t} \in M_{n,1}$.

\begin{defn}\label{defA3}

$M^{<}(B,t,\psi) = \{\ve \in B |\exists ~p~ \in \Z,\q \in \Z^{n}$
such that
   \[\left|
\begin{array}{ll}
    |p + \q\cdot(\ve,\vex  A)| < \psi(2^t) ; \\
    |q_i + \al_iq_n| < 1 ~\forall ~1 \leq i \leq n-1 ; \\
    2^t \leq \| \q \| < 2^{t + 1}. \\
\end{array}%
\right\}. \]
\end{defn}
In a completely analogous fashion, we can define :
\begin{defn}\label{defA4}
$M^{\geq}(B,t,\psi) =  \{\ve \in B |\exists ~p~ \in \Z,\q \in
\Z^{n}$ such that \[\left|
\begin{array}{ll}
    |p + \q \cdot(\ve,\vex  A)| < \psi(2^t) ; \\
    |q_i + \al_iq_n| \geq 1 ~for~ some ~~1 \leq i \leq n-1 ; \\
    2^t \leq \| \q \| < 2^{t + 1}. \\
\end{array}%
\right\}.    \]
\end{defn}
By the Borel-Cantelli Lemma, it is enough to demonstrate :
\begin{claim}\label{cl1}
$\sum_{t = 0}^{\infty}|M^{\geq}(B,t,\psi)| < \infty$.
\end{claim}
and
\begin{claim}\label{cl2}  $\sum_{t =
0}^{\infty}|M^{<}(B,t,\psi)| < \infty$.
\end{claim}

Claim \ref{cl1} turns out to be straightforward, and we provide a
proof in the next section. Claim \ref{cl2} however is more
involved, as one would expect and we devote the rest of the paper
to establishing it.

\section {Establishing Claim \ref{cl1} }
To prove Claim \ref{cl1}, we will establish :\\
\begin{thm}\label{thmA4}
$|M^{\geq}(B,t,\psi)| < C(n,B) \psi(2^{t})2^{nt} $.\\
where $C(n,B)$ is a constant depending on $n$ and $B$.\\
\end{thm}
That the above theorem implies \ref{cl1} is immediate because
\[\sum_{t = 0}^{\infty}|M^{\geq}(B,t,\psi)| \leq \sum_{t =
0}^{\infty} \psi(2^{t})2^{nt} \asymp \int_{0}^{\infty} x^{n -
1}\psi(x)dx \asymp \sum_{k = 0}^{\infty}k^{n - 1}\psi(k) < \infty.
\]
To prove theorem \ref{thmA4}, we will prove :\\
\begin{lemma}\label{thmA5}
For fixed $\q \in \Z^n$
satisfying:\\
\begin{enumerate}
\item $|q_i + \al_iq_n| \geq 1$ for some $i$ and \\
\item $2^t \leq \|\q\| < 2^{t + 1} $.\\
\end{enumerate}
 and for an
arbitrary (fixed) $\theta > 0$, the set:\\
$\hat{M}(B,\q,\theta) =  \{\ve \in B ~|~|p +\q \cdot(\ve,\vex
\cdot A) | < \theta $~
for some $p \in \Z \}$\\
has measure at most $C(n,B) \theta$.\\
\end{lemma}

Clearly  Lemma \ref{thmA5} implies Theorem \ref{thmA4} because\\
$|M^{\geq}(B,t,\psi)| = \sum_{\q \in \Z^n,~2^t \leq \|\q\| < 2^{t
+ 1}} |\hat{M}(B,\q,\psi(2^{t}))| \leq 2^{n +
1}C(n,B) \psi(2^{t})2^{nt} $.\\
And so, it remains to provide the following
\begin{proof}
Let $S = \sqrt{|q_1 + \al_1q_n|^2 + \dots + |q_{n - 1} + \al_{n -
1}q_n|^2}$.~Then $\hat{M}(B,\q,\theta)$ is the union of ``strips"
each of which have volume at most $ C \frac{\theta}{S}(diam B)^{n
- 1}$~ for some constant $C$, depending on $n$.~Also, the number
of such ``strips" is at
most $ S \cdot diam (B) + 1 $.\\
Therefore, $|\hat{M}(B,\q,\theta)|$ is at most $C
\frac{\theta}{S}(diam B)^{n - 1} \times (S \cdot (diam B) + 1)$.\\
From Definition \ref{defA4}, we have $S^{-1} \leq 1 $ and this
implies that
\[|\hat{M}(B,\q,\theta)|< \theta((diam B)^n + (diam B)^{n
- 1}) < \theta C(n,B) \] for a constant $C(n,B)$  as claimed.
\end{proof}
This concludes the easier half of the proof. ~ Note that this
method is completely general and the \hy\  can in fact be replaced
by any affine subspace. A similar situation arises while dealing
with non-degenerate manifolds as well. See \cite{BKM} for details.

\section {Quantitative non-divergence}
To demonstrate Claim \ref{cl2}, we use a variation of a method
used in \cite{BKM} and originally due to D.Kleinbock and
G.A.Margulis, (see \cite{KM}). The method consists of translating
the problem of estimating the measure of $M^{<}(B,t,\psi)$ to an
estimate on non-divergence of lattice-trajectories in Euclidean
space. This section is devoted to setting notation and stating a
theorem which reduces Claim \ref{cl2} to an estimate from
\cite{BKM}.
\begin{defn}
Let $C$ and $\al$ be positive numbers and $V$ be a subset of $\R^d$.
A function $f : V \to \R$ is said to be $(C,\al)-good$ on $V$ if for
any open ball $B \subseteq V$,~and
for any $\ep > 0$, one has :\\
\[\bigg | \bigg\{ \ve \in B \big| |f(\ve)| < \ep \cdot
\sup_{\ve \in B}\{|f(\ve)|\} \bigg \} \bigg | \leq C\ep^{\al}|B|.\]
\end{defn}
Some easy properties of $(C,\al)-good$ functions are :\\
\begin{enumerate}
\item $f$ is $(C,\al)-good$ on $V \Rightarrow$ so is $\lambda f
~\forall ~\lambda ~\in ~\R$.\\
\item $f_i ~~i \in I$ are $(C,\al)-good$ $\Rightarrow$ so is
$\sup_{i \in I}|f_i|$.\\
\end{enumerate}
This paper deals with \hy s, which are parametrized by linear
functions. Consequently, the following lemma will be useful.
\begin{lemma}\label{lem1}
Let $V \subset \R^{d}$, and $f$ be a (continuous) linear function
on $V$. Then $f$ is $(C_{d},1)-good $ on $V$, where $C_{d} =
\frac{2^{d + 2}}{v_{d}} $, and $v_d$ is the volume of the unit
ball in $\R^{d} $.
\end{lemma}
\begin{proof}
Let $B \subseteq V $  be a $d$ dimensional cube,i.e. a product of
$d$ open intervals $B_i$, all of the same size.Let $b$ denote the
side-length of $B$. The set $\{\ve \in B \big| |f(\ve)| < \ep\}$ is
a strip of width at most $2
\frac{\ep b}{\sup_{\ve \in B}|f(\ve) |}  $.\\
The volume of the strip therefore, is at most $2 \frac{\ep
b}{\sup_{\ve \in B}|f(\ve)|}\cdot 2 b^{d - 1}  = 4 \frac{\ep
b^{d}}{\sup_{\ve \in B}|f(\ve)|} = 4 \frac{\ep |B|}{\sup_{\ve \in B}
|f(\ve) |} $ . To complete the proof, we circumscribe open balls
with cubes which are still contained in $V$.
\end{proof}
For more examples of $(C,\al)-good$ functions, see \cite{KM}.\\
Let $\La$ be a discrete subgroup of $\R^k$. A subgroup $\Ga$ of
$\La$ is said to be primitive in $\La$ if $\Ga = \Ga_{\R} \cap
\La$. Let $\Le(\La)$ be
the set of all nonzero primitive subgroups of $\La$.\\
We need to ``measure" discrete subgroups of $\R^l$. Let $\Ga
\subseteq \R^l$ be one such subgroup.~Denote by $\Ga_{\R}$, the
minimal linear subspace of $\R^l$ containing $\Ga$.~If $k =
\dim(\Ga_{\R})$, for any basis $\vv_1, \dots, \vv_k$ of $\Ga$ the
vector $\vv_1 \wedge \dots \wedge \vv_k \in \bigwedge^{k}(\R^l)$ is
(up to a sign) independent of the basis. It is therefore natural to
define  $\|\Ga\| = \|\vw\|$ where $\vw = \vv_1 \wedge \dots
\wedge \vv_k$ is said to represent $\Gamma$. \\
The norm above is the extended norm. More precisely, for $I = (i_1,
\dots,i_j) \subset \{1,\dots,l\},~i_1 < i_2 < \dots < i_j $,~let
\begin{center}
$e_{I} = e_{i_1}\wedge \dots \wedge e_{i_j}
~\in~\bigwedge^{j}(\R^{l})$
\end{center}

Then, for \[\vw = \sum_{I \subset \{1,\dots,l \}} \vw_I e_I,\] we
define $\|\vw\| = $ max$_{I \subset \{1,\dots,k \}}|\vw_I|~$ with
the additional convention that $e_{\emptyset} = 1$.\\
Let $\psi_{0}(x) = x ^{-n}$. Then, from the monotonicity of $\psi$
and the condition in Theorem \ref{thmA3}, it follows that there
exists $x_0 \in  \R $ such that $x > x_0 \Rightarrow \psi(x) <
\psi_{0}(x) $. Hence, $M^{<}(B,t,\psi) \subseteq M^{<}(B,t,\psi_{0})
$, for sufficiently large $t$. Therefore, Claim \ref{cl2} will
follow if we are able to obtain an estimate on $M^{<}(B,t,\psi_{0})
$. And so we can henceforth assume without any loss of generality,
that $\psi(x) = \psi_{0}(x). $ \\ For $0 < \ep < 1$,~let \[D(\ep,t)
= diag(\frac{\ep}{2^{-nt}},
\ep,\dots,\ep,\frac{\ep}{2^{t}},\dots,\frac{\ep}{2^{t}})\]
 i.e. $ D(\ep,t)$ is an element of $\GL(2n,\R$)and has $n$ factors of the type
$\frac{\ep}{2^{t}}$. Thus, $D(\ep,t)$ when applied to a vector
expands the first component and contracts the remaining.\\
In what follows, $\widetilde{A}\in M_{n-1,1}$ will denote the
matrix $(\al_1,\dots,\al_{n - 1})^{t}$. Also, $I_k$ will stand for
the identity $k \times k$ square matrix, and $0$ will stand for
the zero matrix
of appropriate dimension.\\
Let $ \un \in \GL(2n,\R)$ denote the (unipotent) matrix :
\[ \un = \begin{pmatrix}
 1 & 0 & \ve & \vex A\\
 0 & I_{n - 1} & I_{n - 1} & \widetilde{A} \\
 0 & 0 & I_{n-1}& 0\\
 0 & 0& 0 & 1
\end{pmatrix}\]

\begin{defn}\label{defA5}
Let $\La$ denote the subgroup of $\Z^{2n}$ consisting of vectors
of the form
\[\left\{\begin{pmatrix}
p\\
0\\
\vdots\\
0\\
\q^{t}
\end{pmatrix} \mid p \in \Z,\q \in \Z^{n} \right\}\]
\end{defn}
Let $\mathbf{\lambda} \in \La$. Then $\un \mathbf{\lambda}$ is the
vector whose components appear in the definition of
$M^{<}(B,t,\psi)$. Therefore, existence of a non-zero $\q \in \Z^n
\backslash \{0\}$ satisfying the conditions enumerated in Definition
\ref{defA3}  would imply the existence of a nonzero element of
$u_{\ve} \mathbf{\lambda}$ in some parallelepiped in $\R^{2n}$.
Using $D(\ep,t)$ we can transform this parallelepiped into a cube
and then use an estimate from \cite{BKM} to find the measure of $\ve
\in B$ for which this happens. Specifically, it can be seen that
$M^{<}(B,t,\psi)$ is a subset of
\begin{center}
$\{\ve \in B | \|D(\ep,t) \un \mathbf{\lambda}\| < \ep $~for
some~$\mathbf{\lambda} \in \La \backslash \{0\}\}.$
\end{center}

We now state a slightly modified version of Theorem $6.2$~ from
\cite{BKM} which provides the essential :

\begin{thm}\label{thmK3}\cite{BKM}
For arbitrary $l \geq 2$, let $\La$ be a discrete subgroup of
$\R^l$ of rank $k$. Further, let a ball $B(\mathbf{x}_0,r_0)
\subset \R^d$ and a continuous map $H : \tilde{B} \to \GL(W)$ be
given, where $\tilde{B}$ stands for $B(\mathbf{x}_0,3^k
r_0)$.~Take $C,\al > 0,~0 < \rho \leq \frac{1}{k}$ and let
$\|\cdot\|$ be the supremum norm on $\bigwedge(\R^l)$ (as defined
before).
~Assume that for any $\Ga \in \Le(\La)$,\\
\begin{enumerate}
\item the function $\ve \to \|H(\ve)\Ga\|$ is $C,\al-good$ on
$\tilde{B}$~and\\
\item $\exists~ \ve \in B$ such that $\|H(\ve)\Ga\| \geq
\rho$.\\
\end{enumerate}
Then for every positive $\ep \leq \rho$ one has :
\[|\{\ve \in B | \|H(\ve)\mathbf{\lambda}\| <
\ep ~~ for ~some~~ \mathbf{\lambda} \in \La \backslash \{0\}\}| <
k(3^d N_d)^k \cdot C (\frac{\ep}{\rho})^{\al}|B|.\]
\end{thm}

Remarks: \begin{enumerate} \item  The fact that unipotent
trajectories do not diverge dates back to the  work of Margulis in
the $70$'s. The above sharp quantitative estimate is from \cite{BKM}
and is a modification of an earlier estimate from \cite{KM}.\\
\item  In \cite{BKM}, the authors state and prove the theorem for the
norm $\|~\| $ replaced by sub-multiplicative functions on exterior
algebra's of arbitrary finite-dimensional vector spaces. For our
purposes, the sup-norm on $\bigwedge(\R^l)$ suffices and this
changes the restrictions on $\rho$ slightly, from $\rho \leq 1$ in
\cite{BKM} to $\rho \leq \frac{1}{k}$ above. The proof is identical
and rather than producing a verbatim repetition, we direct the
interested reader to \cite{BKM} and to \cite{KSS} for more results
in a similar vein and their applications.
\end{enumerate}

At this stage, we use the fact that we have a Diophantine condition
on the coefficients of the hyperplane, namely, (\ref{defA2}) which
tells us that for some $\delta > 0 $,
\begin{equation}\label{defA14}
 max_{i}|p_i + \al_{i}q  | > |q|^{-n + \delta} ~
\end{equation}
for every $\p \in \Z^{n} $, and all but finitely many $q \in \Z$.
We now state a Theorem which reduces Claim \ref{cl2} to Theorem \ref{thmK3} :\\

\begin{thm}\label{thmA6}
For every $\ep
> 0$ and  for every $t > 0$,\\
\[inf_{\Ga \in \Le(\La)}sup_{\ve \in B}\|D(\ep,t)
\un \Ga\| \geq \ep \cdot min(C_{B}^{1}2^{\frac{\delta t}{n + 1 -
\delta}}, 2^{(n - 1)t}\ep^{n}, C_{B}^{3}2^{t}\ep^{n - 1} )\]\\
where $D(\ep,t)$ and $\un$ are as defined before, $\delta $ is as
in (\ref{defA14}), $\La $ is as in Definition \ref{defA5} and
$C_{B}^{i} ~~ i = 1,3 $ are constants depending on $B$ only.
\end{thm}
To see that Theorem \ref{thmK3} and Theorem \ref{thmA6} imply Claim
\ref{cl2} consider the following
\begin{proof}
We use Theorem \ref{thmK3} with $H(\ve) = D(\ep,t) \un$, and $\ep =
\ep(t) = 2^{-\beta t}$ for some $\beta \in \R^{+} $. It now suffices
to check that the conditions necessary to apply Theorem \ref{thmK3}
are satisfied. From Lemma \ref{lem1}, it follows that condition $1$
is satisfied with $(C,\al) = (C_{n - 1},1)$. As for condition $2$,
notice that a choice of $\beta <$~min$(\frac{n -1}{n +
1},\frac{\delta}{n + 1 - \delta},\frac{1}{n})$~ensures that
\[ \ep \cdot min(C_{B}^{1}2^{\frac{\delta t}{n
+ 1 - \delta}},2^{(n - 1)t}\ep^{n},C_{B}^{3}2^{t}\ep^{n - 1} ) \geq
1 \geq \frac{1}{k}. \] Thus, it follows that :
\[ |M^{<}(B,t,\psi)| < \tilde{C}~\frac{\ep}{\rho}~|B|\]
where $\rho = \frac{1}{k} $ and $\tilde{C} = C_{n - 1}k(3^{n -
1}N_{n - 1})^{k} $ is a constant depending on $n$ and $k$ only (for
details on the constant $N_{n-1}$, see \cite{KM} ). And so,
\[\sum_{t = 0}^{\infty}|M^{<}(B,t,\psi)| < \sum_{t =
0}^{\infty}~\frac{\ep}{\rho} = \sum_{t = 0}^{\infty}
max((C_{B}^{1})^{-1}2^{\frac{-\delta t}{n + 1 - \delta}},2^{-(n -
1)t}\ep^{-n},(C_{B}^{3})^{-1}2^{-t}\ep^{-(n - 1)} )
\]which always converges. Claim \ref{cl2} follows.
\end{proof}
\section {Theorem \ref{thmA6}}
It turns out that obtaining an estimate on $sup_{\ve \in
B}\|D(\ep,t)\un \vw\|$ needs different analysis depending on the
rank of the discrete
subgroup $\vw$ represents.\\
In what follows, we will denote the standard basis of
$\R^{2n}$ by\\
$\{\ee_{0},\ee_{*1},\dots,\ee_{*(n-1)},\ee_1,\dots,\ee_n\}$,
and  $\R^{n + 1} $ will denote the subspace :\\
$\{(x_0,0,\dots,0,x_1,\dots,x_n) ~|~ x_i \in \R,~0 \leq x_i \leq n
\} = \R\ee_0 + \R\ee_1
+ \dots + \R\ee_{n}$.\\
We divide the analysis into the following cases :
\begin{enumerate}
\item $\vw \in \bigwedge^{1}(\R^{n + 1})$.\\
\item $\vw \in \bigwedge^{n + 1}(\R^{n + 1}).$\\
\item $\vw \in \bigwedge^{k}(\R^{n + 1})$ for $1 <k \leq n$.
\end{enumerate}

The action of $\un$ on these basis vectors is as follows :
\begin{itemize}
\item $\un\ee_0 = \ee_0.$\\
\item $\un\ee_{*i} = \ee_{*i}.$\\
\item $\un\ee_i =  x_i\ee_0 + \ee_{*i} + \ee_i,~~1 \leq i \leq n - 1.$\\
\item $\un\ee_n =  (\al_0 + \sum_{i = 1}^{n - 1} \al_ix_i)\ee_0 +
\sum_{i = 1}^{n - 1}\al_i
\ee_{*(i)} + \ee_n.$\\
\end{itemize}
Case $1.$ Let $\vw \in \bigwedge^{1}(\R^{n + 1})$
be of the form \\
$\vw = p_0\ee_0 + q_1\ee_1 + \dots + q_{n - 1}\ee_{n - 1}
+ q_n\ee_n$. Then,\\
$D(\ep,t)\un \vw = \frac{\ep}{2^{-nt}}(p_0 + q_1x_1 + \dots + q_{n
- 1}x_{n - 1} + q_n(\al_0 + \al_1x_1 + \dots + \al_{n - 1}x_{n -
1}))\ee_0 + \ep(q_1 + q_n\al_1)\ee_{*1} + \dots + \ep(q_{n - 1} +
q_n\al_{n - 1})\ee_{*(n - 1)} + \frac{\ep}{2^{t}}q_1\ee_1 + \dots
+ \frac{\ep}{2^{t}}q_n\ee_n$.\\
Recall that we have a \di\ condition on the coefficients of the
\hy,~namely, (\ref{defA2}) which tells us that for some $\delta
> 0 $,
\[ max_{i}|p_i + \al_{i}q  | > |q|^{-n + \delta} ~ \] for every $\p \in \Z^{n}
$, and all but finitely many $q \in \Z$. This readily implies that
for all but finitely many $\q \in \Z^{n} $, we have
\begin{equation} \label{defA6}
 max_{i} |q_i + \al_{i}q_n | > |q_n|^{-n + \delta}
\end{equation}
The coefficient of $\ee_0$ in $D(\ep,t)\un w$ can be written as
\[\frac{\ep}{2^{-nt}}(p_0 + q_n\al_0 + (q_1 +
\al_1q_n)x_1 + \dots + (q_{n-1} + \al_{n-1}q_n)x_{n-1}).\] From
(\ref{defA6}), we can then deduce that for every $\ve \in B$,we have
\begin{equation}\label{defA7}
|\frac{\ep}{2^{-nt}}(p_0 + q_n\al_0 + (q_1 + \al_1q_n)x_1 + \dots
+ (q_{n-1} + \al_{n-1}q_n)x_{n-1})| >
\frac{\ep}{2^{-nt}}C_{B}|q_n|^{-n + \delta}
\end{equation}
where $C_B$ is a constant depending on the ball $B$.\\
We now turn our attention to the terms $\frac{\ep}{2^{t}}q_i\ee_i ~1
\leq  i \leq n-1 $. Our strategy is to choose a non-zero $ q_i $ and
then compare the coefficient of $\ee_i$, with the term containing
$q_i$ in the coefficient of $\ee_0$. Due to (\ref{defA7}), it makes
sense to instead consider the equation:
\begin{equation}\label{defA8}
\frac{\ep}{2^{-nt}}C_{B}y^{-n + \delta} = \frac{\ep}{2^{t}}y.
\end{equation}
We will denote the unique root of (\ref{defA8}) by $y_0$. So $y_0$
is in fact, equal to $(C_B)^{\frac{1}{n + 1 - \delta}}[2^{t(n +
1)}]^{\frac{1}{ n + 1 - \delta}} $. From (\ref{defA7}) and
(\ref{defA8}),it follows that for $w \in \bigwedge^{1}(\R^{n +
1})$,
\begin{equation}\label{defA9}
 sup_{\ve \in B}\|D(\ep,t)\un \vw \| \geq
C_{B}^{1}\ep \cdot 2^{t(\frac{\delta}{n + 1 - \delta})}
\end{equation}
where $C_{B}^{1} $ is a constant depending on $B$ only.\\
Let us now consider Case $2$ i.e.\\
$\vw = \ee_0 \wedge \ee_1 \wedge \dots \wedge \ee_n
~\in~\bigwedge^{n + 1}(\R^{n + 1})$. In this case, $\un \vw $
contains the term $\ee_0 \wedge \ee_{*1} \wedge \ee_{*2} \wedge
\dots \wedge
\ee_{*(n - 1)}\wedge \ee_n$. Therefore,\\
$sup_{\ve \in B}\|D(\ep,t) \un \vw\| \geq 2^{(n - 1)t}\cdot \ep^{n
+ 1}$, which implies that for $\vw \in ~\bigwedge^{n + 1}(\R^{n +
1})$,
\begin{equation}\label{defA10}
  sup_{\ve \in B}\|D(\ep,t) \un \vw\| >
\ep \cdot 2^{(n - 1)t} \cdot \ep^{n}.
\end{equation}
Now for Case $3$. Here, it is enough to project to the subspace of
$\R^{2n}$ spanned by $\ee_0,\ee_1,\dots,\ee_n$. To simplify
computations, we instead consider $\R^{n + 1}$ and (with some abuse
of notation), denote its standard basis by
$\ee_0,\ee_1,\dots,\ee_n$. For this, we must modify our setup a
little. Let $\hat{u}_{\ve}$ denote the matrix :
\[\hat{u}_{\ve} = \begin{pmatrix}
 1 & \ve & \vex A\\
 0& I_{n-1} & 0 \\
 0 & 0 & 1
\end{pmatrix}\]
In what follows, it will be convenient to change notation slightly.
Define functions $\hat{f}_i(\ve)$ for $\ve \in B$, in the following
manner:
\[ \hat{f}_i(\ve) = x_i   ~~1 \leq i < n.\]
\[ \hat{f}_n(\ve) = \vex A.\]
Further, let $\mathbf{\hat{f}}(\ve) =
(\hat{f}_{1}(\ve),\dots,\hat{f}_{n}(\ve)).$ Then, $\hat{u}_{\ve} $
can now be written as:
\[\hat{u}_{\ve} = \begin{pmatrix}
 1 & \mathbf{\hat{f}}(\ve)\\
 \tilde{0}&  I_n \\
\end{pmatrix} \]
Let $\hat{D}(\ep,t) =
diag(\frac{\ep}{2^{-nt}},\frac{\ep}{2^{t}},\dots,\frac{\ep}{2^{t}})$
where, as before, $0 < \ep < 1 $. Thus, both $\hat{u}_{\ve}$ and
$\hat{D}(\ep,t)$ are matrices in $\GL(n + 1,\R)$.
It therefore suffices to obtain an estimate on $ \\
sup_{\ve \in B}\|\hat{D}(\ep,t)\hat{u}_{\ve}\vw\|$ for any $\vw
\in \bigwedge^{j}(\R^{n + 1})$.\\
Let us write down an expression for $\hat{u}_{\ve} \ee_I$.\\
$\hat{u}_{\ve}$ acts on the basis vectors as follows :\\
\begin{enumerate}
\item $\hat{u}_{\ve} \ee_0 = \ee_0.$\\
\item $\hat{u}_{\ve} \ee_i = \hat{f}_{i}(\ve)\ee_0 + \ee_i ~ 1
\leq i \leq n$.
\end{enumerate}

Therefore, for $\ee_{I} = \ee_{i_1}\wedge \dots \wedge \ee_{i_j}
~\in~\bigwedge^{j}(\R^{n + 1})$, we have :
\begin{itemize}
\item $\hat{u}_{\ve} \ee_I = \ee_I $~if $0 \in I$.\\
\item $\hat{u}_{\ve} \ee_I = \ee_I + \sum_{i \in I}
(-1)^{l(I,i)}\hat{f}_i(\ve) \ee_{I \cup \{0 \} \backslash \{i\}}$
\end{itemize}

where $l(I,i) =$ number of elements of $I$ strictly between $0$ and
$i$.
Taking $\vw$ of the form $\vw = \sum_{I}\vw_I \ee_I$ we get\\
$\hat{u}_{\ve} \vw = \sum_{0 \in I}(\vw_I + \sum_{i \notin
I}(-1)^{l(I,i)} \vw_{I \cup \{i\} \backslash
\{0\}}\hat{f}_{i}(\ve))\ee_I + \sum_{0
\notin I}\vw_I\ee_I$.\\
For each $\vw = \sum_{I}\vw_I \ee_I$, following \cite{K}, we define
a vector $\mathbf{c}_{I,\vw} \in \R^{n + 1}$ as follows:
\[\mathbf{c}_{I,\vw} = \sum_{i \notin (I \backslash \{0\})}
(-1)^{l(I,i)} \vw_{I \cup \{i\} \backslash \{0\}}.\] Then, the
non-constant components of $\hat{u}_{\ve} \vw $ can
be written as $(1,\hat{f}_1,\hat{f}_2,\dots,\hat{f}_n)\mathbf{c}_{I,\vw} $.\\
Let $P$ be the augmented $(n,n + 1)$ matrix given by $(I_{n}|A)$.
Then, we can write \[(1,\hat{f}_1,\dots,\hat{f}_n) =
\mathbf{\tilde{x}}P.\] Therefore, the nonconstant components of
$\hat{u}_{\ve} \vw $ can now be
written as $\mathbf{\tilde{x}}P\mathbf{c}_{I,\vw} $.\\
We can thus replace $\|\hat{D}(\ep,t)\hat{u}_{\ve} \vw \| $ by
\begin{equation}\label{defA11}
C^{3}_{B}max(\ep^{j}2^{t(n + 1 - j)} max_{0 \in
I}\|P\mathbf{c}_{I,\vw} \|,\ep^{j}2^{-jt}max_{0 \notin I}|\vw_I| )
\end{equation}
where $C^{3}_{B}$ is a constant depending on $B$ alone.\\
At this stage, it is advantageous to split the vector
$\mathbf{c}_{I,\vw} $ into components, $\mathbf{c}_{I,\vw}^{+} \in
\R^{n}$ and $\mathbf{c}_{I,\vw}^{-} \in \R $ (which respectively
denote first $n$ rows of the column vector $\mathbf{c}_{I,\vw} $ and
its last row). And so we can write :
 \[\mathbf{c}_{I,\vw} = \begin{pmatrix}
 \mathbf{c}_{I,\vw}^{+}\\
 \mathbf{c}_{I,\vw}^{-} \\
\end{pmatrix}\]
where $\mathbf{c}_{I,\vw}^{+} = \vw_I \ee_0 + \sum_{i \in
\{1,\dots,n - 1\} \backslash I} (-1)^{l(I,i)} \vw_{I \cup \{i\}
\backslash \{0\}}
\ee_i$~and\\
$\mathbf{c}_{I,\vw}^{-} = (-1)^{l(I,n)} \vw_{I \cup \{n\}
\backslash \{0\}}$ if
$n \notin I$ and $0$ otherwise.\\
Therefore, $P\mathbf{c}_{I,\vw} = (\mathbf{c}_{I,\vw}^{+} + A
\mathbf{c}_{I,\vw}^{-}) $, and it follows that :
\begin{equation}\label{defA12}
sup_{\ve \in B}\|\hat{D}(\ep,t)\hat{u_{\ve}}\vw \| \geq C_{B}^{3}
\{\ep^{j} 2^{t(n + 1 - j)} max_{0 \in I}\|\mathbf{c}_{I,\vw}^{+} +
A \mathbf{c}_{I,\vw}^{-} \|,\ep^{j}2^{-j t} max_{0 \notin
I}|\vw_I| \}
\end{equation}
where $C_{B}^{3}$ is as in (\ref{defA11}). The stage is now set for
Lemma $4.5$ from \cite{K} which tells us that:
\begin{lemma}\cite{K}\label{thmK4}
For any $\vw \in \bigwedge^{j}(\R^{n + 1})$, we have :
\[max_{ 0 \in I } \|\mathbf{c}_{I,\vw}^{+} + A \mathbf{c}_{I,\vw}^{-}\| \geq 1.\]
\end{lemma}
From (\ref{defA11}), (\ref{defA12}) and Lemma \ref{thmK4}, it
follows that:
\begin{equation}\label{defA13}
 sup_{\ve \in B}\|\hat{D}(\ep,t) \hat{u}_{\ve} \vw\| \geq
C_{B}^{3} 2^{t(n - j + 1)}\cdot \ep^{j} \geq C_{B}^{3} \cdot \ep
\cdot 2^{t} \cdot \ep^{n - 1}.
\end{equation}

From (\ref{defA9}), (\ref{defA10}) and (\ref{defA13}), it follows
that for any $\Ga \in \Le(\La) $, we have
\[ \|D(\ep,t)
u_{\ve}\Ga\| \geq \ep \cdot min(C_{B}^{1}2^{\frac{\delta t}{n + 1 -
\delta}},2^{(n - 1)t}\ep^{n },C_{B}^{3}2^{t}\ep^{n - 1}).
\] This completes the proof of Theorem \ref{thmA6} and hence the
proof of Claim \ref{cl2}.$\qed$

\section { Conclusion and Open Questions}

As we have mentioned, in \cite{BKM}, the authors have obtained a
Khintchine type theorem for non-degenerate manifolds of any
codimension. Moreover, \cite{BBDD} treats the case of a straight
line passing through the origin. We could thus hope to complete
the picture and prove a theorem similar to Theorem \ref{thmA3} for
affine subspaces of arbitrary codimension. Following \cite{K}, it
would also be nice to establish a Khintchine theorem for
non-degenerate submanifolds of affine subspaces. This will be
dealt with in a forthcoming paper, where we also establish the
multiplicative analogue of the convergence \kh\ theorem.\\
It would also be of interest to check if the Diophantine condition
given in this paper is optimal. More conjectures along these lines
can be found in \cite{K}.\\\\{\bf Acknowledgements.} The author
thanks his advisor Dmitry Kleinbock for suggesting the problem, for
constant encouragement and for numerous helpful discussions.

\address{Anish Ghosh\\MS $050$, Brandeis
University\\Waltham, MA-$02454$\\U.S.A.}\\
\email{ghosh@brandeis.edu}

\end{document}